\documentclass[final]{siamltex}

\usepackage{epsfig}
\usepackage{amsmath}
\usepackage{amssymb}
\usepackage{subfigure}

\newtheorem{algorithm}{Algorithm}

\def\remark{{\sc Remark:\hspace{0.5em}}}

\def\eps{\varepsilon}

\def\cond{\operatorname{cond}}
\def\diag{\operatorname{diag}}
\def\ord{{\cal O}}
\def\toep{\operatorname{toep}}
\def\Zst{{\mathbb Z}}

\def\Nst{{\mathbb N}}

\def\Rst{{\mathbb R}}

\def\Psp{{\boldsymbol P}}
\def\PM{{\Psp_{\! M}}}
\def\PMM{{\Psp_{\! M_x M_y}}}

\def\kron{{\otimes}}

\def\Taug{{T_{\text{aug}}}}
\def\Haug{{H_{\text{aug}}}}
\def\Aaug{{A_{\text{aug}}}}
\def\xaug{{x_{\text{aug}}}}
\def\yaug{{y_{\text{aug}}}}
\def\sw{{s^{(w)}}}
\def\nullmat{{{\bf 0}}}
\def\Lsp{{\boldsymbol L}}
\def\Ltsp{{\Lsp^2}}

\def\lsp{{\boldsymbol\ell}}
\def\ltsp{{\lsp^2}}

\def\shadowbox{\hbox{\rule[-0.0ex]{0.1ex}{1.2ex}%
\hspace{-0.1ex}\rule[-0.0ex]{1.2ex}{0.1ex}%
\hspace{0.0ex}\rule[-0.0ex]{0.1ex}{1.2ex}\hspace{-1.3ex}%
\rule[1.15ex]{1.25ex}{0.1ex}\hspace{-0.0ex}\rule[-0.25ex]{0.3ex}{1.1ex}%
\hspace{-1.2ex}\rule[-0.25ex]{1.1ex}{0.25ex}}}
\def\QED{\mbox{\phantom{m}}\nolinebreak\hfill$\,\shadowbox$}
\def\proof{{\it Proof. }}

\def\adots{\mathinner{\mkern1mu\raise1pt\hbox{.}\mkern2mu
   \raise4pt\hbox{.}\mkern2mu\raise7pt\vbox{\kern7pt\hbox{.}}\mkern1mu}}

\begin{document}

\title{\bf Fast multi-dimensional scattered data approximation with
Neumann boundary conditions}

\thanks{The second author has been supported in part by NSF grant 9973373.}

\author{Denis Grishin and Thomas Strohmer\thanks{Department of Mathematics, 
        University of California, Davis, CA 95616-8633, USA;
        Correspondence: strohmer@math.ucdavis.edu.}}
\date{}
\maketitle

\vspace*{-5cm}
\noindent
%{\footnotesize \textsc{submitted to SIAM J.\ Sci.\ Comp.} \\
{\footnotesize \textsc{} \\

\vspace*{4.5cm}

\begin{abstract}
An important problem in applications is the approximation of a
function $f$ from a finite set of randomly scattered data $f(x_j)$. 
A common and powerful approach is to construct a trigonometric least squares
approximation based on the set of exponentials $\{e^{2\pi i kx}\}$. This 
leads to fast numerical algorithms, 
but suffers from disturbing boundary effects due to the underlying 
periodicity assumption on the data, an assumption that is rarely satisfied 
in practice. To overcome this drawback we impose Neumann boundary
conditions on the data. This implies the use of cosine polynomials
$\cos (\pi kx)$ as basis functions. We show that scattered data
approximation using cosine polynomials leads to a least squares 
problem involving certain Toeplitz+Hankel matrices. We derive estimates on 
the condition number of these matrices.  Unlike other 
Toeplitz+Hankel matrices, the Toeplitz+Hankel matrices arising in our
context cannot be diagonalized by the discrete cosine 
transform, but they still allow a fast matrix-vector multiplication via
DCT which gives rise to fast conjugate gradient type algorithms. 
We show how the results can be generalized to higher
dimensions. Finally we demonstrate the performance of the proposed
method by applying it to a two-dimensional geophysical scattered data problem.
\end{abstract}


\noindent

\pagestyle{myheadings}
\thispagestyle{plain}
%\markboth{DENIS GRISHIN AND THOMAS STROHMER}{MULTI-DIMENSIONAL SCATTERED
%DATA APPROXIMATION}

\section{Introduction}
\label{intro}
An ubiquitous problem in mathematics and in applications is the
reconstruction or approximation of a function $f$ from its non-uniformly 
spaced sampling values $s_j = f(x_j)$. 
Without further knowledge about $f$ this is an ill-posed problem, since
the subspace of functions $h$ with $h(x_j)=s_j$ has always
infinite dimension. 
Moreover in practice we are given only a finite
number of samples $\{s_j\}_{j=1}^r$, which makes a complete reconstruction
of $f$ in general impossible, so the best we can hope for is to
compute a good approximation to $f$.
Fortunately in many practical situations the functions under consideration
are not arbitrary, but possess some smoothness properties.
For instance physics often implies that $f$ is bandlimited.
In this and many other cases a linear combination of 
trigonometric basis functions $\{e^{2\pi i kx}\}_{k \in \Zst}$ often provides 
a good approximation to $f$. 
Other powerful models for scattered data approximation are based
on radial basis functions and on shift-invraint systems~\cite{Sch95}.
%The exponentials $e^{2\pi kx}$ can also be intepreted as finite-dimensional
%model for the infinite-dimensional space of bandlimited functions. 

Least squares approximation using exponentials as basis functions provides
a tool that is general enough to be useful in a variety of situations where 
smooth functions are involved, while the algebraic structure of the functions
$e^{2\pi ikx}$ is rich enough to give rise to fast and robust numerical 
algorithms to compute the approximation, cf.~e.g.~\cite{RAG91,FGS95,Fas97}.

Arguably the main drawback of approximation by exponentials is the
underlying periodicity assumption about the function to be approximated.
To be more precise, let $f$ be a smooth continuous function and let 
$\{f(x_j)\}_{j=1}^r$ be samples of $f$ taken at the points
$x_1 < \dots < x_r$. Without loss of generality we assume that $x_1=0$ and
$x_r=1$. We want to approximate $f$ on the sampling interval 
$[x_1,x_r) = [0,1)$ by a trigonometric polynomial 
$p(x)=\sum_{k=-M}^{M} c_k e^{2\pi ikx}$ with $M < r/2$. If $f(0)=f(1)$ we 
can safely conclude from Weierstrass' theorem that a trigonometric polynomial
of low degree will give a good approximation to $f$ on the interval
$[0,1)$. However if $f(0) \neq f(1)$ then this difference is felt
as discontinuity by the approximating polynomial $p$. In fact standard 
Fourier analysis tells us that the coefficients $\{c_k\}_{k \in \Zst}$ of 
$p$ will at best decay like {\it o}$(1/k)$, thus a large degree $M$ is
required to obtain a reasonable approximation to $f$ on $[0,1)$.
However since in practice only a finite number of samples is available
we may not be able to choose $M$ sufficiently large to obtain a
satisfactory approximation to $f$.

A standard method to enforce periodicity of $f$ on $[0,1)$ is to 
multiply $f$ with a smooth ``window function'' $w$ which decays
rapidly to zero at the boundaries of the sampling interval. However
such a procedure can considerably reduce the interval in which the
approximation is in agreement with the ``non-windowed'' sampling values
$f(x_j)$. We could also try to reduce the unpleasant behavior caused by 
the boundary effects
by choosing the period $N$ of $p$ slightly larger than the length of the
sampling interval. Nevertheless, if $|f(0)-f(1)|$ is large we still
need a polynomial of large degree to obtain a reasonable approximation
to $f$ on $[0,1)$. We also note that boundary effects become worse 
with increasing dimension.

Instead of extending $f$ (respectively its samples $f(x_j)$) periodically
across the boundaries of the sampling interval, we can apply Neumann boundary
conditions, i.e., a symmetric extension across the end points of the
sampling interval. This has the big advantage that we avoid the discontinuity 
at the boundaries. The Fourier coefficients of a continuous (periodic) 
function decay at least like {\it o}$1/k$ and at best like {\it o}$1/k^2$. 
Thus loosely spoken, the decay is one order of magnitude faster than
compared to a periodic extension.
This faster decay implies that a lower polynomial degree should suffice
to obtain a good trigonometric approximation.\footnote{This is exactly the
reason why the (old) JPEG image compression algorithm uses the DCT instead
of the DFT.}

If we extend the sampling values $f(x_j)_{j=1}^r$
symmetrically across the boundaries we obtain a sampling sequence
that is periodic on the interval $[0,2)$ and symmetric with
respect to the midpoint $1$. To adapt the trigonometric basis
functions to this situation we have to replace the exponentials
$\{e^{2\pi i kx}\}_{k \in \Zst}$ by the basis functions 
$\{\cos(\pi kx)\}_{k \in \Nst}$. The functions $\cos(\pi kx)$
are symmetric around 1 and periodic with respect to the interval $[0,2)$.
The advantage when using cosine polynomials instead of exponentials
is obvious from the discussion above: we reduce disturbing boundary
effects, which results in a better approximation of the original function.

In the case of trigonometric approximation based on exponentials
it has been shown that the least squares approximation can be formulated
as hermitian positive definite Toeplitz system~\cite{FGS95}. 
Gr\"ochenig has derived explicit
bounds for the condition number of the Toeplitz matrix 
that allow to estimate the stability and convergence
of the involved numerical algorithms~\cite{Gro93a,FGS95}.
Moreover all steps to compute and solve the Toeplitz system can 
be done quickly by (nonuniform) FFT-based methods.

The crucial questions that we will investigate in this paper are: 
Does the least squares approximation problem using cosine polynomials
also give rise to a linear system of equation whose matrix has a nice
structure? Can we find fast and robust numerical algorithms to solve the 
least squares problem? Can we give a priori estimates on the condition number
of the matrix? Can we generalize the algorithm easily to higher dimensions?
How does our approach perform for real world problems?
This paper is devoted to clarify these questions. 

The rest of the paper is organized as follows. In Section~\ref{s:two} we 
analyze the least squares approximation problem using cosine polynomials. We 
show that the resulting matrix has a certain Toeplitz+Hankel structure and
derive estimates on the condition number of this matrix. 
In Section~\ref{s:three} we present a fast algorithm to solve the least squares
problem using the conjugate gradient method and the discrete cosine
transform (DCT). The generalization to the multi-dimensional case is 
described in 
Section~\ref{s:four}. Finally in Section~\ref{s:five} we demonstrate the 
performance of the proposed method by applying it to a scattered data 
problem arising in geophysics.

%We illustrate this by a simple example. Let $f(x)=\sum_{k=-M}^{M}
%a_k e^{2\pi i kx}$. We have sampled $f$ at the randomly spaced points 

\medskip

The idea of using Neumann boundary conditions instead of periodic boundary 
conditions has turned out to be very fruitful in the context of image 
deblurring problems. In fact, the research presented in this paper was 
inspired by the article {\em A fast algorithm for deblurring models with 
Neumann boundary conditions} by Michael Ng, Raymond Chan, and 
W.C.~Tang~\cite{NCT99}.

\section{Nonuniform sampling, cosine polynomials, and Toeplitz+Hankel matrices} 
%\section{Fast scattered data approximation using cosine polynomials} 
\label{s:two}

%\subsection{Scattered data approximation and Toeplitz+Hankel matrices} 
%\label{ss:th}

We start by defining the space $\PM$ of cosine polynomials of maximal 
degree $M$ as
\begin{equation}
\PM =\left\{p: p(x) = \frac{c_0}{\sqrt{2}}+\sum_{k=1}^{M} c_k \cos(\pi kx),
 c = \{c_k\}_{k=0}^{M} \in \Rst^{M+1} \right\}.
\label{cospol}
\end{equation}
There are two reasons for the introduction of the $1/\sqrt{2}$-scaling factor 
of the coefficient $c_0$ in~\eqref{cospol}. The first reason is that
we have the Parseval type identity 
\begin{equation}
\|p\|^2_2 = \int \limits_{-\infty}^{+\infty}|p(x)|^2 dx =
\frac{c_0^2}{2} +\frac{1}{2} \sum_{k=1}^{M} c_k^2 = \frac{1}{2} \|c\|^2_2.
\label{parseval}
\end{equation}
The second reason is increased stability of the numerical
algorithms we are going to derive, as we will 
explain in the remark after Theorem~\ref{th:cond}.

Let us return to the approximation problem. Given sampling 
points\footnote{Throughout the paper we will always assume that
the sampling locations $x_j$ are pairwise distinct.}
$\{x_j\}_{j=1}^{r}$ and sampling values $\{s_j\}_{j=1}^{r}$, we want to 
solve the least squares problem
\begin{equation}
\underset{p \in \PM}{\min} \sum_{j=1}^{r}|p(x_j) - s_j|^2 w_j.
\label{lsp1}
\end{equation}
Here the $w_j>0$ are weights which the user may choose at her
convenience. Often the trivial choice $w_j=1$ is sufficient. In other cases
it is useful to choose the weights such that they compensate for
irregularities in the sampling set, i.e., smaller weights are used
in regions with high sampling density and larger weights in regions with
few sampling points. In~\eqref{lsp1} we have assumed that the polynomial
degree $M$ is fixed. We will discuss the important question of
how to determine the appropriate degree of the approximating polynomial
in Section~\ref{s:three}.

By defining the $r \times (M+1)$ Vandermonde-like matrix $V$ via
\begin{equation}
\label{vander}
V_{j,k} = 
\begin{cases}
\frac{1}{\sqrt{2}} \sqrt{w_j}, & \text{for $k=0; j=1,\dots,r$},\\
 \sqrt{w_j}\cos (\pi k x_j),   & \text{for $k=1,\dots M; j=1,\dots,r$},
\end{cases}
\end{equation}
and setting $\sw = \{\sqrt{w_j} s_j\}_{j=1}^{r}$ 
we can reformulate the least squares problem~\eqref{lsp1} as
\begin{equation}
\underset{c \in \Rst^{M+1}}{\min} \|Vc - \sw\|^2_2.
\label{lsp2}
\end{equation}

It is well-known that the solution of \eqref{lsp2} can be computed
by solving the normal equations
\begin{equation}
\label{normal}
V^{T} V c = V^{T} \sw.
\end{equation}
Switching to the normal equations can lead to problems
of numerical instability due to the squaring of the condition number
of $V$. However, as we will see, the system matrix
of the normal equations has a very nice algebraic structure that paves the
way to fast numerical algorithms for solving~\eqref{lsp1}.
Thus to handle the trade-off between numerical
stability and computational efficiency it is important to have
an a priori estimate of the condition number of the matrix $V$.
Such an estimate will aid us in the decision if we shall compute
the least squares solution by a direct solution of the system $Vc=\sw$ or by
switching to the system $V^{T}V c = V^{T}\sw$.

The following theorem provides both insight in the algebraic
structure of $V^{T} V$ and an upper bound of
the condition number of $V^{T} V$.

\begin{theorem}
\label{th:cond}
Assume we are given nonuniformly spaced sampling points 
$\{x_j\}_{j=1}^r \in [0,1]$, sampling
values $s=\{s_j\}_{j=1}^r$ and positive weights $\{w_j\}_{j=1}^r$.
Define $A:= V^{T} V$, where $V$ is as in~\eqref{vander},
and set $b=V^{T}\sw$. There holds: \\
(i) The matrix $A$ is a scaled Toeplitz+Hankel matrix of the form
\begin{equation}
A = D (T+H) D, 
\label{dthd}
\end{equation}
where
\begin{equation}
T = 
\begin{bmatrix}
a_0     & a_1     & \dots  & a_{M-1} & a_{M}   \\
a_1     & a_0     & \ddots &         & a_{M-1} \\
\vdots  & \ddots  & \ddots & \ddots  & \vdots  \\
a_{M-1} &         & \ddots & \ddots  & a_{1}   \\
a_M     & a_{M-1} & \dots  & a_1     & a_0
\end{bmatrix},\,\,
H = 
\begin{bmatrix}
a_0     & a_1     & \dots  & a_{M-1} & a_{M}   \\
a_1     & a_2     & \adots & a_M     & a_{M+1} \\
\vdots  & \adots  & \adots & \adots  & \vdots  \\
a_{M-1} & \adots  & \adots &         & a_{2M}  \\
a_M     & a_{M+1} & \dots  & a_{2M}  & a_{2M+1}
\end{bmatrix},
\label{toephank}
\end{equation}
with 
\begin{equation}
\label{Aentries}
a_k = \frac{1}{2}\sum_{j=1}^{r} w_j \cos(\pi k x_j), 
\qquad k=0,\dots,2M+1,
\end{equation}
and $D=\diag (\frac{1}{\sqrt{2}},1,\dots,1)$. \\
(ii) If $M < r$ then $A$ is invertible and the coefficient vector 
$c=\{c_k\}_{k=0}^{M}$ of the cosine polynomial $p \in \PM$ 
that solves~\eqref{lsp1} is given by
\begin{equation}
c = A^{-1}b.
\label{Ainv}
\end{equation}
(iii) Define the weights $w_j$ by
\begin{equation}
w_j = \frac{x_{j+1}-x_{j-1}}{2}, \qquad j=1,\dots,r,
\label{weights}
\end{equation}
where we set $x_{0}:=-x_1, x_{r+1}:=2-x_r$. If
\begin{equation}
\label{maxgap}
\delta := \max_j |x_{j+1} - x_{j}| < \frac{1}{M}
\end{equation}
then the condition number $\kappa (A)$ is bounded by
\begin{equation}
\label{Acond}
\kappa (A) \le \frac{(1+\delta M)^2}{(1-\delta M)^2}.
\end{equation}
\end{theorem}

%\begin{proof}
\proof
(i) Note that 
\begin{equation}
\label{p1}
A_{k,l} = (V^{T} V)_{k,l}
= \eps_{k,l}\sum_{j=1}^{r} w_j \cos (\pi l x_j) \cos (\pi k x_j),
\qquad k,l=0,\dots,M,
\end{equation}
where 
\begin{equation}
\label{p2}
\eps_{k,l}=
\begin{cases}
\frac{1}{2} & \text{if $k=0$ and $l=0$}, \\
\frac{1}{\sqrt{2}} & \text{if $k=0$ or $l=0$, $k \neq l$}, \\
1 & \text{if $k>0$ and $l>0$.} 
\end{cases}
\end{equation}
The result follows now readily from a simple calculation by applying 
the formula
\begin{equation}
\cos (\alpha) \cos(\beta) = \cos(\alpha + \beta) + \cos(\alpha-\beta),
\label{p3}
\end{equation}
to~\eqref{p1} and using the fact that the entries of $T$ and $H$ satisfy 
$T_{k,l}= a_{k-l}$ and $H_{k,l} = a_{k+l}$ respectively.

(ii) The invertibility of $A$ follows from the well-known fact that
the Vandermonde-like matrix $V$ has rank $M+1$ for mutually
different points $x_j$ (assuming $w_j \neq 0$).
The rest follows from~\eqref{normal}.

(iii) With the exception of a few minor modifications the proof of this
part is similar to Gr\"ochenig's elegant proof on the upper bound of the 
condition number of certain Toeplitz matrices, see~\cite{Gro93a}. However
instead of confronting the reader with a patchwork of required 
modifications of Gr\"ochenig's proof we prefer to present a complete proof.

The proof makes use of {\em Wirtinger's inequality}~\cite{HLP52}:
If $f \in \Ltsp(a,b)$ and either $f(a)=0$ or $f(b)=0$, then
\begin{equation}
\int \limits_{a}^{b} \big|f(x)\big|^2 \, dx \le \frac{4}{\pi^2}(b-a)^2
\int \limits_{a}^{b} \big|f'(x)\big|^2 \, dx .
\label{wirt}
\end{equation}

We proceed with the proof of (iii).
Let $P$ be the orthogonal projection of $\Ltsp([0,1])$ onto
$P_M$. Define the operator $S$ by
\begin{equation}
Sp = P \Big( \sum_{j=1}^{r} p(x_j) \chi_j \Big).
\label{cond1}
\end{equation}
Here $\chi_j(x)$ denotes the characteristic function of the
interval $[y_{j-1}, y_j]$, where $y_j = \frac{x_{j+1}-x_j}{2}, 
j=1,\dots, r$ with $x_0=-x_1, x_{r+1}=1-x_r$.

We compute
\begin{gather}
\|p - Sp\|^2_2 = \|P\big( \sum_{j=1}^{r}(p - p(x_j) ) \chi_j \big)\|^2_2
\le \| \sum_{j=1}^{r}(p - p(x_j) ) \chi_j \|^2_2 = \notag\\
 \int \limits_{0}^{1}\big|\sum_{j=1}^{r}(p - p(x_j) ) \chi_j \big|^2 dx
= \sum_{j=1}^{r}\, \int \limits_{y_{j-1}}^{y_j}|p - p(x_j)|^2 dx.
\label{cond2}
\end{gather}
We write
$$
\int \limits_{y_{j-1}}^{y_j}|p - p(x_j)|^2 dx =
\int \limits_{y_{j-1}}^{x_j}|p - p(x_j)|^2 dx +
\int \limits_{x_j}^{y_j}|p - p(x_j)|^2 dx,
$$
and apply Wirtinger's inequality~\eqref{wirt} to each of the integrals on 
the left-hand side. Since $|y_j - x_j| \le \delta/2$ and 
$|x_j - y_{j-1}| \le \delta/2$ we obtain
\begin{equation}
\sum_{j=1}^{r}\int \limits_{y_{j-1}}^{y_j}|p - p(x_j)| dx \le
\frac{\delta^2}{\pi^2} \sum_{j=1}^{r}\int \limits_{0}^{1}|p'(x)| dx 
= \frac{\delta^2}{\pi^2} \|p'\|^2_2.
\label{cond3}
\end{equation}
Note that
\begin{equation}
p'(x) = \sum_{k=0}^{M} c_k \pi k \sin(\pi k x),
= \sum_{k=1}^{M} c_k \pi k \sin(\pi k x).
\label{cond4}
\end{equation}
Hence we have the Bernstein type inequality
\begin{gather}
\|p'\|^2_2=\int \limits_{0}^{1} |\sum_{k=1}^{M} c_k \pi k \sin(\pi k x)|^2 dx
\notag \\
\le (\pi M)^2 \int \limits_{0}^{1} |\sum_{k=1}^{M} c_k \sin(\pi k x)|^2 dx
\le (\pi M)^2 \|p\|^2_2. 
\label{cond5}
\end{gather}
Thus by combining \eqref{cond2},~\eqref{cond3} and \eqref{cond5} we get
\begin{equation}
\|p -Sp\|^2_2 \le \delta^2 M^2 \|p\|^2_2.
\label{cond6}
\end{equation}
Hence
\begin{equation}
\|I-S\|_{\text{op}} \le \delta M,
\label{cond7}
\end{equation}
and since $\delta < 1/M$ by assumption, we conclude that $S$ is invertible
and
\begin{equation}
\|S^{-1}\|_{\text{op}} \le (1-\delta M)^{-1}.
\label{cond8}
\end{equation}
There holds
\begin{gather}
(1-\delta M)^2 \|p\|^2_2 = (1-\delta M)^2 \|S^{-1}Sp\|^2_{\text{op}} \le 
\notag \\
\le (1-\delta M)^2 \|S^{-1}\|^2_{\text{op}} \|Sp\|^2_2 \le \|Sp\|^2_2 
\le \sum_{j=1}^{r} |p(x_j)|^2 w_j.
\label{cond9}
\end{gather}
Also
\begin{gather}
\sum_{j=1}^{r} |p(x_j)|^2 w_j \le \|p-p+\sum_{j=1}^{r} p(x_j) \chi_j \|^2_2
\le \big(\|p\| + \|p-\sum_{j=1}^{r} p(x_j) \chi_j \|_2\big)^2 \notag \\
\le \big(\|p\|_2 + \delta M \|p\|_2 \big)^2 \le (1+\delta M)^2 \|p\|^2_2. 
\label{cond10}
\end{gather}
Thus
\begin{equation}
(1-\delta M)^2 \|p\|^2_2 \le \sum_{j=1}^{r} |p(x_j)|^2 w_j
\le (1+\delta M)^2 \|p\|^2_2.
\label{cond11}
\end{equation}
By definition we have for any $p \in \PM$ with coefficient vector $a$
\begin{equation}
\langle Aa,a \rangle = \langle V^{T}Va,a \rangle =
\langle V a,V a\rangle = \sum_{j=1}^{r} |p(x_j)|^2 w_j.
\label{cond12}
\end{equation}
Using the relation $\|p\|^2_2 = \frac{1}{2}\|a\|^2_2$ we obtain
\begin{equation}
\frac{1}{2}(1-\delta M)^2 \|a\|^2_2 \le \langle Aa,a \rangle
\le \frac{1}{2}(1+\delta M)^2 \|a\|^2_2,
\label{cond13}
\end{equation}
and therefore
$$\kappa (A) \le \frac{(1+\delta M)^2}{(1-\delta M)^2}.$$
\QED
%\end{proof}

\remark 
%In order to make transparent why we introduced the $1/\sqrt{2}$
%scaling of $c_0$ in the definition of cosine polynomials, 
We briefly analyze the least squares problem~\eqref{lsp2} when using non-scaled
cosine polynomials $\tilde{p}(x) = \sum_{k=0}^{M} c_k \cos(\pi kx)$.
It is easy to see that the corresponding Vandermonde-like matrix $\tilde{V}$
satisfies 
\begin{equation}
\tilde{V}D = V,
\label{r1}
\end{equation}
with $D$ as in part~(i) of Theorem~\ref{th:cond} and $V$ as
in~\eqref{vander}. Hence 
\begin{equation}
\tilde{A} := \tilde{V}^{T}\tilde{V} = D^{-1} V^{T}VD^{-1}. 
\label{r2}
\end{equation}
The estimates
\begin{equation}
\|\tilde{A}x\|_2 \le \|D^{-1}\|^2_{\text{op}} \|A\|_{\text{op}} \|x\|_2 \le 2
\|A\|_{\text{op}} \|x\|_2,
\label{r3}
\end{equation}
and
\begin{equation}
\|\tilde{A}^{-1}x\|_2 \le \|D\|^2_{\text{op}} \|A\|_{\text{op}} \|x\|_2 \le
\|A\|_{\text{op}} \|x\|_2,
\label{r4}
\end{equation}
imply that
\begin{equation}
\cond(\tilde{A}) \le 2 \cond(A).
\label{r5}
\end{equation}
Thus the condition number of $\tilde{A}$ can be twice as large as
the condition number of $A$. This is why we prefer to use scaled cosine
polynomials as defined in~\eqref{cospol}. The inequality~\eqref{r5} is sharp 
as can be seen from the following simple example. Let the sampling points 
$x_j$ be equally spaced, and choose the weights $w_j$ as in
Theorem~\ref{th:cond}. In this case it is not difficult to see that
\begin{equation}
A = \frac{1}{2} I_{M+1},
\label{r6}
\end{equation}
where $I_{M+1}$ denotes the $(M+1) \times (M+1)$ identity matrix, whereas
\begin{equation}
\tilde{A} = \frac{1}{2}
\begin{bmatrix}
2      & 0 & \dots   & 0      \\
0      & 1 & \dots   & 0      \\
\vdots &   & \ddots  & \vdots \\
0      &   &         & 1
\end{bmatrix}.
\label{r7}
\end{equation}
Thus obviously $\cond(\tilde{A}) = 2\cond(A)$ in this case.

%\section{Fast solution of the least squares problem} \label{ss:three}
\section{Fast computation of the least squares approximation} \label{s:three}

In this section we present a fast algorithm for solving the least squares
problem~\eqref{lsp1}. Our algorithm is based on the conjugate
gradient method in connection with a fast matrix-vector multiplication
involving the DCT.
Before we proceed we briefly review some properties of the DCT-I.
There are four types of the DCT, cf.~\cite{Str99d}. For our purposes
we will use the (scaled) DCT-I.

\begin{definition}
\label{defdct1}
The {\em Type-I Discrete Cosine Transform matrix} 
(DCT-I for short) of size $n \times n$ is defined by
\begin{equation}
[C_n]_{k,l} = 
\begin{cases}
\frac{1}{\sqrt{2n-2}}\cos (\pi \frac{kl}{n-1}) & \text{if $k=0$ or $k=n-1$}, \\
\frac{2}{\sqrt{2n-2}}\cos (\pi \frac{kl}{n-1}) & \text{if $k=1,\dots,n-2$}.
\end{cases}
\label{dct1}
\end{equation}
If the dimension of the matrix $C_n$ is clear from the context we
drop the subscript and simply write $C$ instead.
\end{definition}

The DCT-I matrix $C$ satisfies $C C = I$. It is not unitary, but
can be easily made unitary by appropriate scaling. For define
the diagonal matrix $\tilde{D} = \diag([1 ,\sqrt{2}, \dots ,\sqrt{2}, 1])$ and
set $\tilde{C} = \tilde{D}{-1} C \tilde{D}$. Then it is easy to see that 
$\tilde{C} \tilde{C}^T = I$. In some cases it is more convenient to work
with $\tilde{C}$ instead of $C$~\cite{HR98}. However the results presented 
in this paper can be more elegantly expressed when using the 
definition~\eqref{dct1} of the DCT-I. Fast algorithms for computing
$Cx$ require 2.5 $\ord (n \log n)$ operations if $x$ is a vector
of length $n+1$ and $n$ is a power of two~\cite{Loa92}, cf.~also
\cite{Ste92,BT97}.

It is well-known that the DCT-I matrix diagonalizes certain Toeplitz+Hankel
matrices~\cite{SPP95,HR98}. For let $T = \toep(a)$ be a symmetric Toeplitz 
matrix with first column $a=[a_0, a_1, \dots, a_n]^T$. We
define the counter-identity matrix $J$ by
\begin{equation}
%T = 
%\begin{bmatrix}
%b_0     & b_1     & \dots  & b_{n-1}   \\
%b_1     & b_0     & \ddots & \vdots    \\
%\vdots  & \ddots  & \ddots & \vdots  \\
%b_{n-1} &         & \dots  & b_0
%\end{bmatrix},\,\,
%%H = 
%\begin{bmatrix}
%b_0     & b_1     & \dots  & b_{n-2} & b_{n-1} \\
%b_1     & b_2     & \adots & b_{n-1} & b_{n-2} \\
%\vdots  & \adots  & \adots & \adots  & \vdots  \\
%b_{n-2} & \adots  & \adots &         & b_1     \\
%b_{n-1} & b_{n-2} & \dots  & b_1     & b_0
%\end{bmatrix}
J = 
\begin{bmatrix}
0 &        & 1  \\
  & \adots &    \\
1 &        & 0
\end{bmatrix}.
\label{counterid}
\end{equation}
If
\begin{equation}
B=\toep(a) +J \toep( Ja) := T+H
\label{th}
\end{equation}
(note that $J \toep(Ja)$ is a Hankel matrix that is symmetric with respect to
the counter diagonal) then
\begin{equation}
\label{dctdiag}
C^T B C = \Sigma, \qquad
\text{where $\Sigma$ is a diagonal matrix.}
\end{equation}

An important consequence of this diagonalization property is that the 
multiplication of a matrix $B$ of 
the form~\eqref{th} with a vector $x$ can be carried out in $\ord (n \log n)$
operations via DCT-I \cite{BT97}, similar to the
multiplication of a vector by a Toeplitz matrix which can be computed via FFT
by embedding the Toeplitz matrix into a circulant matrix.

To be precise, assume we want to 
compute $y=Bx$ where $C^{T}BC = \Sigma$. There holds
\begin{equation}
y = Bx = C^T C^TBCCx = C^T \Sigma Cx.
\label{ybx}
\end{equation}
Of course in a numerical implementation we would not compute the
diagonal matrix $\Sigma$ explicitly. Instead we proceed as follows. Let 
$b$ be the first column of $B$, define the scaling 
matrix $D_1=\diag(2,1,\dots,1,2)$ and observe that $C = D_1^{-1}C^T  D_1$.
A simple calculation shows that $D_1^{-1}\Sigma = 
\sqrt{\frac{n-1}{2}}\diag(C^T b)$.
Hence 
$$y=Bx=C^T D_1 D_1^{-1} \Sigma Cx = \sqrt{\frac{n-1}{2}} 
C^T \diag(D_1 C^T b ) D_1^{-1}C^T  D_1 x, $$
and therefore
\begin{equation}
y=\sqrt{\frac{n-1}{2}}C^T \big[ (C^T b)\circ (C^T D_1 x])\big],
\label{ybx1}
\end{equation}
where the operation ``$\circ$'' denotes the pointwise product
between vectors. Hence the product $Bx$ can be computed by three DCT-I's
in $\ord (n \log n)$ operations.

Observe that the Toeplitz+Hankel part of the matrix $A=D(T+H)D$ in~\eqref{dthd} 
of Theorem~\ref{th:cond} is {\em not} of the form~\eqref{th}, since the first 
row and the last column of the Hankel matrix $H$ in~\eqref{toephank} have 
different entries. Thus $A$ is {\em not} diagonalized by the DCT-I (or any
other DCT). But we can embed the Toeplitz+Hankel part of $A$ in a 
Toeplitz+Hankel matrix of the form~\eqref{th}, similar to the embedding of 
a Toeplitz matrix in a circulant matrix. To see this, let $T$ and $H$ be 
defined as in~\eqref{toephank}. We embed $T+H$ in the $(2M+1) \times
(2M+1)$ augmented Toeplitz+Hankel matrix $\Taug+\Haug$, where
\begin{gather}
\Taug = 
\begin{bmatrix}
a_0 & \dots & a_{M} & a_{M+1} & \dots & a_{2M}  \\
\vdots & \ddots & \ddots & \ddots & \ddots & \vdots  \\
a_{M} & \ddots & \ddots & \ddots & \ddots & a_{M} \\
a_{M+1} &\ddots & \ddots & \ddots & \ddots & a_{M-1} \\
\vdots & \ddots & \ddots & \ddots & \ddots & \vdots  \\
a_{2M} & \dots & a_{M+1} & a_M & \dots & a_0
\end{bmatrix},\\
\Haug = 
\begin{bmatrix}
a_0 & \dots & a_{M} & a_{M+1} & \dots & a_{2M}  \\
\vdots & \adots & \adots & \adots & \adots & \vdots  \\
a_{M} & \adots & \adots & \adots & \adots & a_{M} \\
a_{M+1} &\adots & \adots & \adots & \adots & a_{M-1} \\
\vdots & \adots & \adots & \adots & \adots & \vdots  \\
a_{2M} & \dots & a_{M+1} & a_M & \dots & a_0
\end{bmatrix}.
\label{}
\end{gather}
The matrix $T+H$ is the $(M+1)\times (M+1)$ principal leading submatrix
of $\Taug+\Haug$.

Thus for a DCT-I based fast implementation of the matrix vector
product $Ax$ we proceed as follows. We write $y=Ax=D(T+H)Dx$ and
define $\xaug := [(Dx)^{T}, 0,\dots, 0]^T$. Compute $\yaug = \Aaug \xaug$
according to~\eqref{ybx1}. The vector $y$ is then given by the first $M+1$ 
entries of $\yaug$ multiplied by $D$.

In order to obtain augmented matrices whose size is $2^n +1$
we can always insert as many zeros as necessary after
$a_{2M}$ in the first row of $\Taug$ and $\Haug$ 
without destroying the algebraic structure of the matrices.
Thus the matrix vector multiplication $Ax$ can always be carried out in 
$\ord (M \log M)$. This zero-padding is similar to the zero-padding
of the Toeplitz case (where the zeros are added in the middle
of the first row).

Note that a direct computation of the entries of the matrix $A$ and of
the right hand side $b$
will take $\ord (Mr)$ operations. Thus, although we can solve the
system $Ax=b$ in $\ord (M \log M)$ operations, the computation of the
entries of $A$ and $b$ will soon become the bottleneck for large scale
problems. Fortunately there exist fast algorithms for computing sums
of the form~\eqref{Aentries}. In~\cite{Pot01} Daniel Potts has developed
fast algorithms for computing the DCT for nonuniformly spaced points.
Like nonuniform FFT algorithms~\cite{PST01} a nonuniform DCT-I
(NDCT for short) can be computed in
$\ord(\alpha M \log(\alpha M) +mr)$ operations, where $\alpha$
and $m$ are constants. See~\cite{Pot01} for details.

Based on the observations above, we propose the following
fast algorithm for solving the least squares problem~\eqref{lsp1}.

\begin{algorithm}[Fast scattered data approximation using cosine polynomials]
\label{alg1}
\\
{\bf Input:} Nonuniformly spaced sampling points 
$\{x_j\}_{j=1}^r \in [0,1]$, sampling values
$\{s_j\}_{j=1}^r$, weights $\{w_j\}_{j=1}^r$ and user-defined points 
$\{t_l\}_{l=0}^L \in [0,1]$.  \\
{\bf Task:} Compute the coefficients of the cosine 
polynomial of degree $M$ that solves~\eqref{lsp1} and evaluate the 
polynomial at the points $\{t_l\}_{l=1}^L$.

\noindent
{\bf Step 1:} Compute the first column of $A$ in~\eqref{dthd} and the
right hand side $b=V^T \sw$ via NDCT. This takes
$\ord(\alpha M \log(\alpha M) +mr)$ operations, where $\alpha$
and $m$ are (small) constants.

\noindent
{\bf Step 2:} Solve $Ac=b$ iteratively by the conjugate gradient method. 
Using fast matrix-vector multiplication this can be done in
$\ord (M \log M)$ operations per iteration.

\noindent
{\bf Step 3:} Evaluate $p(x)=\frac{c_0}{\sqrt{2}}+\sum_{k=1}^{M} 
c_k \cos (\pi k x)$ at the points $\{t_l\}_{l=0}^L$. If $t_l = l/L$
and $L=2^n$ for some $n \in \Nst$, then this can be done by a DCT
in $\ord (L \log L)$ operations. If $L \neq 2^n$ we can use a fast
radix-$p$ DCT, see~\cite{Ste92}. If the $t_l$ are nonuniformly
spaced we use a NDCT to compute $\{p(t_l)\}_{l=1}^L$.

\noindent
{\bf Output:} Least squares approximating polynomial $p$ of degree $M$, 
evaluated at the points $\{t_l\}_{l=1}^L$.
\end{algorithm}

\remark
If the sampling set satisfies the maximal gap condition~\eqref{maxgap}
and the weights are chosen according to~\eqref{weights} we can utilize
the bound on $\kappa(A)$ in~\eqref{Acond} of Theorem~\ref{th:cond} to 
estimate the rate of CG using the standard formula~\cite{GL96}
\begin{equation}
\|c^{(n)} - c\|_2 \le 2 \kappa(A) \left(
\frac{\sqrt{\kappa(A)}-1}{\sqrt{\kappa(A)}+1}\right)^k \|c^{(0)} - c\|_2,
\label{cgest}
\end{equation}
where $c^{(n)}$ denotes the solution after the $n$-th iteration of CG
applied to $Ac=b$.

If the condition number of $A$ is large (whether or not the maximal
gap condition is satisfied) it may be better to solve the least squares 
problem~\eqref{lsp1} $Vc=b$ without explicitly establishing the normal 
equations. One can resort to ``non-symmetric'' versions of CG such
as GMRES or LSQR, cf.~\cite{GL96}. Since the NDCT provides a fast way to 
carry out the multiplication of the matrix $V$ with a vector we still obtain
a fast algorithm. However the computational costs are in general
larger than those for Algorithm~\ref{alg1} since a NDCT is more expensive
than a DCT and the NDCT has to be applied in each iteration, whereas in
Algorithm~\ref{alg1} it has to be applied only in the initial stage of
the algorithm.

If the matrix $A$ is ill-conditioned due to large gaps in the sampling set
one might be tempted to apply one of the cosine-transform based 
preconditioners to improve the situation. However preconditioners cannot
significantly improve the stability in this case. This can be shown in a 
similar way as it is done in Section~4.2 of~\cite{Str99c} for 
trigonometric approximation using exponentials.

There exist fast direct methods to solve Toeplitz+Hankel systems
(not all of them apply to our situation though), 
see~\cite{KS99} and in particular the work of Heinig~\cite{HR98,Hei01} . 
But many of these solvers require that the matrix dimension is a power of two.
It is possible to overcome this severe constraint, however at the
cost of a more involved algorithm. As we have seen for the 
conjugate gradient iterations the initial size of the matrix does
not play a major role, since when constructing the augmented matrix we can
always insert the appropriate number of zeros to get a size
of a power of two. Furthermore, if the set of sampling points is a jittered
version of a set of regularly spaced points, standard perturbation
theory implies that the eigenvalues of $A$ will be clustered around
1. Thus CG will converge in very few iterations. Direct solvers
cannot take advantake of such sitations.

\subsection{Multilevel scattered data approximation} \label{ss:}

The reader may have noticed that we have tacitly assumed that the
polynomial degree $M$ is given a priori. Although this is a common
assumption in polynomial approximation it is not justified in 
many applications. In fact, the appropriate choice of $M$ has a major
influence on the usefulness of the resulting approximating
polynomial, cf.~\cite{Str98}.
In~\cite{SS98} Otmar Scherzer and the second author have developed
a multilevel scheme that automatically adapts to the solution 
of the optimal ``level'', i.e., the optimal polynomial degree in our case.
This multilevel algorithm applies to our approximation method without
modification. 

In a nutshell the multilevel version of Algorithm~\ref{alg1}
works as follows, for details we refer to~\cite{SS98,GS01}. We start at 
the first {\it level} with an initial choice for the approximating 
polynomial (e.g., $M_0=1$) and apply Algorithm~\ref{alg1}.
We stop the CG iterations when a specific stopping criterion is 
satisfied and obtain the approximation $p_1$, say. Then we proceed to the 
next {\it level} by choosing a degree $M_1 > M_0$ (e.g., $M_1 = M_0+1$). 
We use the approximation $p_1$ from the previous level as initial guess for
the solution at the new level and apply Algorithm~\ref{alg1}. We proceed 
through increasing levels until at the $k$-th level the approximating 
polynomial $p_k$ satisfies the {\it discrepancy principle} 
\begin{equation}
\sum_{j=1}^{r} |p_k(x_j)-s_j| w_j \le \eps \sum_{j=1}^{r} |s_j|^2 w_j,
\label{discrepancy}
\end{equation}
where $\eps$ is a parameter related to the accuracy of the given data
$s_j$.

A fast $\ord (M \log M)$ implementation of the multi-level
scheme for cosine polynomials can be derived in a similar way as it is
done for the exponentials, see Algorithm~2 in Section 5.1 of~\cite{GS01}.
An crucial observation thereby is that the scaled Toeplitz+Hankel matrix 
$A_M$ associated with the least squares problem~\eqref{lsp2} for
degree $M$ is related to the matrix $A_{M+1}$ associated with
the least squares problem~\eqref{lsp2} in a nice way. Namely, 
$A_M$ is the principal leading submatrix of $A_{M+1}$.

\bigskip
\remark Finding the optimal level for the approximating function
is a common and important problem in scattered data approximation.
When using radial basis functions or shift-invariant systems as model
one has to deal with the trade-off between accuracy and stability when
determining the width of the basis functions,
cf. e.g.\cite{Sch95}. The multi-level idea provides a natural
framework to handle this trade-off.

\section{Two-dimensional scattered data approximation} \label{s:four}

Many of the results of the previous sections can be extended to
arbitrary dimensions. For the sake of simplicity of notation we will
focus mainly on the two-dimensional case.

We are given sampling values $s=\{s_j\}_{j=1}^r$ and 
randomly spaced sampling points $\{(x_j,y_j)\}_{j=1}^r$. Without loss
of generality we assume that $(x_j,y_j) \in [0,1] \times [0,1]$, otherwise 
we can always renormalize the sampling points accordingly.

The space $\PMM$ consists of two-dimensional cosine polynomials $p$
of degree $M_x M_y$ defined by
\begin{equation}
p(x,y) = \frac{c_{0,0}}{\sqrt{2}}+
\underset{\max\{k,l\}>0}{\sum_{k=0}^{M_x}\sum_{l=0}^{M_y}}
 c_{k,l} \cos(\pi kx) \cos(\pi ly),
\label{cospol2}
\end{equation}
with real-valued coefficients $c_{k,l}$.

Analogous to the one-dimensional scattered data problem we want to
find the $p \in \PMM$ that solves
\begin{equation}
\min \sum_{j=1}^{r} |p(x_j,y_j) - s_j|^2 w_j.
\label{s1}
\end{equation}
We define the block matrix $V$ by
\begin{gather}
V = 
\begin{bmatrix}
 V^{(0)} & V^{(1)} & \dots & V^{(M_y)}
\end{bmatrix},\label{s2a} \\
\text{with} \,\,V^{(l)}_{j,k} = 
\eps_{k,l}\sqrt{w_j} \cos(\pi kx_j) \cos(\pi l y_j), \quad j=1,\dots,r,
\label{s2b} \\
\text{where $\eps_{k,l}$} =
\begin{cases}
\frac{1}{\sqrt{2}} & \text{if $k=0$ and $l=0$,} \\
1 & \text{if $k=0,\dots,M_x; l=0,\dots,M_y; \max\{k,l\}>0$.} 
\end{cases}
\end{gather}
By stacking the columns of $c$ and with a slight 
abuse of notation we can rewrite~\eqref{s1} as 
\begin{equation}
\min \| Vc - \sw\|,
\label{s3}
\end{equation}
where $\sw = \{\sqrt{w_j} s_j\}_{j=1}^r$.

Similar to the 1-D case, we can solve~\eqref{s3} by switching
to the normal equations. The next theorem describes the algebraic
structure of the system matrix of the normal equations.

\begin{theorem}
\label{th:blockth}
Let $V$ be as defined in~\eqref{s2a}-\eqref{s2b}. Then the matrix 
$A:=V^T V$ is a scaled block Toeplitz+Hankel matrix of the form $A=D(T+H)D$ 
with
\begin{gather}
T =
\begin{bmatrix}
A^{(0)}   & A^{(1)} & \dots  & A^{(M_y-1)} & A^{(M_y)} \\
A^{(1)}   & A^{(0)} & \ddots &           & A^{(M_y-1)} \\
\vdots    & \ddots  & \ddots & \ddots    & \vdots    \\
\vdots    & \ddots  & \ddots & \ddots    & \vdots    \\
A^{(M_y)} & \dots   & \dots  & \dots     & A^{(0)}
\end{bmatrix}, \\
H=
\begin{bmatrix}
A^{(0)}  & A^{(1)} & \dots  & A^{(M_y-1)}  & A^{(M_y)}  \\
A^{(1)}  & A^{(2)} & \adots & \adots     & A^{(M_y+1)}    \\
\vdots   & \adots  & \adots & \adots     & \vdots     \\
\vdots   & \adots  & \adots &            & A^{(2M_y-1)} \\
A^{(M_y)} & \dots  & \dots  & A^{(2M_y-1)} & A^{(2M_y)}
\end{bmatrix},
\label{bthmat}
\end{gather}
where each block $A^{(k)}, k=0,\dots,2M_y$ is an $(M_x+1) \times (M_x+1)$
matrix of the form $A^{(k)}=T^{(k)}+H^{(k)}$ with $T{(k)}$ and $H^{(k)}$
as in~\eqref{toephank} and $D=\diag (\frac{1}{\sqrt{2}},1,\dots,1)$. 
\end{theorem}

%\begin{proof}
\proof
It follows from~\eqref{s2a} and~\eqref{s2b} that 
\begin{gather}
A_{l,l',k,k'} = \eps_{k,l} \eps_{k',l'} \sum_{j=1}^{r} w_j \Big(
\cos (\pi (l+l')x_j) \cos( \pi (k+k') y_j) +
\cos (\pi (l-l')x_j) \cos( \pi (k+k') y_j) + \notag \\
\qquad \qquad + \cos (\pi (l+l')x_j) \cos( \pi (k-k') y_j) +
\cos (\pi (l-l')x_j) \cos( \pi (k-k') y_j)\Big).
\label{entries}
\end{gather}
Here the indices $l,l'$ refer to the $(l,l')$-th block of $A$ and the
indices $k,k'$ refer to the element in the $k$-th row and $k'$-th
column in a certain block.

Now we consider the entries of $A$ for fixed $l$ and $l'$. Using 
formula~\eqref{p3} we calculate
\begin{align}
A_{l,l',k,k'} =  \eps_{k,l} \eps_{k',l'} \sum_{j=1}^{r} w_j \Big(
& c_1 [\cos (\pi (k+k')y_j) + \cos( \pi (k-k') y_j)] + \notag \\
+ & c_2 [\cos (\pi (k+k')y_j) + \cos( \pi (k-k') y_j)]\Big), 
\qquad k,k'=0,\dots,M_x, 
\label{entries2}
\end{align}
where the constants $c_1$ and $c_2$ are given by 
$c_1:=\cos(\pi (l+l')x_j)$, $c_2:=\cos(\pi (l-l')x_j)$.
Thus the $(l,l')$-th block of $A$ is indeed of the form~\eqref{dthd}.

By repeating this step with reversed roles for $k,k'$ and $l,l'$
we see that the ``global'' structure of $A$ is of the form~\eqref{bthmat}.
\QED
%\end{proof}

In order to utilize the block Toeplitz+Hankel structure of the normal 
equations we have to extend the fact that the DCT-I diagonalizes certain
Toeplitz+Hankel matrices to the case of block Toeplitz+Hankel matrices. 

We need some preparation before we proceed.
Let $B$ be a block matrix of the form
\begin{equation}
B =
\begin{bmatrix}
B^{(0,0)} & \dots & B^{(0,n-1)} \\
\vdots    &       & \vdots      \\
B^{(n-1,0)} & \dots & B^{(n-1,n-1)} 
\end{bmatrix}
\label{blockmat}
\end{equation}
where the blocks $B^{(k,l)}$ are matrices of size $m \times m$. 
For such block matrices we
define the {\em mod-$m$ permutation matrix} $\Pi_{m,n}$ via
\begin{equation}
[\Pi_{m,n} B \Pi^T_{m,n}]_{i,j;k,l} = B_{k,l;i,j}, \qquad
0 \le i,j \le m-1, 0 \le k,l \le n-1.
\label{}
\end{equation}
In words, the $(i,j)$-th entry of the $(k,l)$-th block of $B$ is permuted
to the $(k,l)$-th entry of the $(i,j)$-th block. We have $\Pi_{m,n} =
\Pi_{n,m}^T$, see~\cite{Loa92}.

\begin{definition}
\label{defdct2}
The {\em two-dimensional type-I Discrete Cosine Transform} 
of an $m \times n$ signal $x$ is given by 
\begin{align}
(C x)_{i,j} = 
\frac{\eps_{i,j}}{\sqrt{2m-2}\sqrt{2n-2}} \sum_{k=0}^{m-1}\sum_{l=0}^{n-1}
x_{k,l} & \cos \big(\pi \frac{ik}{m-1}\big)\cos \big(\pi \frac{jl}{n-1}\big),
\label{dct2} \\
& \qquad i=0,\dots,m-1; j=0,\dots,n-1,
\end{align}
where
\begin{equation}
\eps_{i,j}=
\begin{cases}
1 & \text{if $i \in \{0,m-1\}$ and $j \in \{0,n-1\}$,} \\
2 & \text{if $i \in \{0,m-1\}$ and $j \notin \{0,n-1\}$,} \\
2 & \text{if $i \notin \{0,m-1\}$ and $j \in \{0,n-1\}$,} \\
4 & \text{if $i =1,\dots,m-2$ and $j=1,\dots,n-2$.}
\end{cases}
\notag
\end{equation}
\end{definition}
The two-dimensional DCT-I can be represented by the $mn \times mn$ matrix
$C_m \kron C_n$ where the matrices $C_m$ and $C_n$ represent one-dimensional
DCT-I's as in definition~\ref{defdct1} and $\kron$ denotes the usual
Kronecker product.

Similar to the 1-D DCT-I the 2-D DCT-I diagonalizes certain
block Toeplitz+Hankel matrices.
\begin{theorem}
\label{th:blockdiag}
A matrix $B$ is diagonalized by a two-dimensional DCT-I if and only if 
$B$ is of the form
\begin{equation}
B =
\begin{bmatrix}
B^{(0)}   & B^{(1)} & \dots  & B^{(n-2)} & B^{(n-1)} \\
B^{(1)}   & B^{(0)} & \ddots &           & B^{(n-2)} \\
\vdots    & \ddots  & \ddots & \ddots    & \vdots    \\
\vdots    & \ddots  & \ddots & \ddots    & \vdots    \\
B^{(n-1)} & \dots   & \dots  & \dots     & B^{(0)}
\end{bmatrix}
+
\begin{bmatrix}
B^{(0)}  & B^{(1)} & \dots  & B^{(n-2)}  & B^{(n-1)}  \\
B^{(1)}  & B^{(2)} & \adots & \adots     & B^{(n-2)}    \\
\vdots   & \adots  & \adots & \adots     & \vdots     \\
\vdots   & \adots  & \adots &            & \vdots \\
B^{(n-1)} & \dots  & \dots  & \dots      & B^{(0)}
\end{bmatrix},
\label{bthdiag}
\end{equation}
where each block $B^{(k)}, k=0,\dots,n-1$ is a $m \times m$ 
Toeplitz+Hankel matrix of the form~\eqref{th}.
\end{theorem}

%\begin{proof}
\proof
The proof is similar to the proof of Theorem~3.3 in~\cite{NCT99} and
uses basic properties of the Kronecker product $\kron$. Let $B$ be a block 
Toeplitz+Hankel matrix as in the assumption of the theorem. We have to show 
that $B$ is diagonalized by the two-dimensional DCT-I $C=C_n \kron C_m$. 
Note that each block $B^{(k)}$ of $B$ can be diagonalized by a one-dimensional 
DCT-I $C_m$, i.e., $C_m^T B^{(k)} C_m =  \Lambda^{(k)}$, $k=0,\dots n-1$, where
the $\Lambda^{(k)}$ are $m \times m$ diagonal matrices. Since 
$C_n \kron C_m = (C_n \kron I_m) (I_n \kron C_m)$ it follows that
\begin{gather}
(C_n \kron C_m)^T B(C_n \kron C_m) = 
(C_n^T \kron I_m) (I_n \kron C_m^T) B (I_n \kron C_m)(C_n \kron I_m)=
(C_n^T \kron I_m) \Lambda (C_n \kron I_m),
\label{}
\end{gather}
where
\begin{equation}
\Lambda  =
\begin{bmatrix}
\Lambda^{(0)}   & \Lambda^{(1)} & \dots  & \Lambda^{(n-2)} & \Lambda^{(n-1)} \\
\Lambda^{(1)}   & \Lambda^{(0)} & \ddots &           & \Lambda^{(n-2)} \\
\vdots    & \ddots  & \ddots & \ddots    & \vdots    \\
\vdots    & \ddots  & \ddots & \ddots    & \vdots    \\
\Lambda^{(n-1)} & \dots   & \dots  & \dots     & \Lambda^{(0)}
\end{bmatrix}
+
\begin{bmatrix}
\Lambda^{(0)}  & \Lambda^{(1)} & \dots  & \Lambda^{(n-2)}  & \Lambda^{(n-1)}  \\
\Lambda^{(1)}  & \Lambda^{(2)} & \adots & \adots     & \Lambda^{(n-2)}    \\
\vdots   & \adots  & \adots & \adots     & \vdots     \\
\vdots   & \adots  & \adots &            & \vdots \\
\Lambda^{(n-1)} & \dots  & \dots  & \dots & \Lambda^{(0)}
\end{bmatrix}.
\label{lambdamat}
\end{equation}
We compute
\begin{equation}
\Pi_{m,n} \Lambda \Pi^T_{m,n} =  \tilde{B} =
\begin{bmatrix}
\tilde{B}^{(0)} & \nullmat        & \dots     & \nullmat          \\
\nullmat        & \tilde{B}^{(1)} &           & \vdots            \\
\vdots          &                 & \ddots    & \nullmat          \\
\nullmat        & \dots           & \nullmat  & \tilde{B}^{(m-1)} \\
\end{bmatrix},
\label{}
\end{equation}
where $\nullmat$ is an $n\times n$ zero matrix.
It follows from~\eqref{lambdamat} that each $\tilde{B}^{(k)}, k=0,\dots,m-1$ 
is an $n \times n$ Toeplitz+Hankel matrix of the form~\eqref{th}. Therefore
$C_n^T \tilde{B}^{(k)} C_n = \tilde{\Lambda}^{(k)},
k=0,\dots,m-1$.

Since $\Pi_{m,n}(C_n^T \kron I_m) \Pi^T_{m,n} = I_m \kron C_n^T$ 
(e.g., see~\cite{Loa92}) we have
\begin{align}
(C_n^T \kron I_m) \Lambda (C_n \kron I_m) = & 
\Pi_{m,n} \Pi^T_{m,n}(C_n^T \kron I_m)\Pi^T_{m,n} \Pi_{m,n} \Lambda
\Pi_{m,n}^T \Pi_{m,n} (C_n \kron I_m)\Pi_{m,n} \Pi^T_{m,n} \notag \\
= & \Pi^T_{m,n} (I_m \kron C_n^T) \tilde{B}
 (I_m \kron C_n^T) \Pi_{m,n} \notag \\
= & \Pi^T_{m,n} \tilde{\Lambda} \Pi_{m,n},
\label{}
\end{align}
where $\tilde{\Lambda}$ is a block diagonal matrix with diagonal
blocks $\tilde{\Lambda}^{(k)}$. Thus $\tilde{\Lambda}$ is a diagonal matrix.
It follows from the definition of $\Pi_{m,n}$ that 
$\Pi^T_{m,n} \tilde{\Lambda} \Pi_{m,n}$ is then also a diagonal
matrix.

The opposite direction follows from the fact that $C C = I$. 
\QED
%\end{proof}

The matrix $A$ associated with the least squares problem
\eqref{s3} is {\em not} diagonalized by the 2-D DCT-I. But analogous
to the 1-D case, $A$ can be embedded into a block Toeplitz+Hankel
matrix that is diagonalized by the 2-D DCT-I. Thus similar to the
1-D case the matrix-vector multiplication $Ax$ can be carried out
in $\ord (M_x M_y \log M_x M_y)$ operations. 

We leave it to the reader to extend Theorems~\ref{th:blockth} and
\ref{th:blockdiag} and the fast matrix-vector multiplication
to dimensions larger than two.
Since the NDCT can also be generalized to two and higher dimensions we have
a fast numerical algorithm for computing the least squares approximation
using cosine polynomials in multiple dimensions in the same
way as it is outlined in Algorithm~\ref{alg1}.

\remark There is one notable difficulty that arises when considering
the scattered data approximation problem in higher dimensions.
In the 1-D case a sufficient condition for invertibility of
the matrix $A$ is that the polynomial degree $M$ is smaller than
the number of samples $r$. This is an immediate consequence of the
fundamental theorem of algebra. Unfortunately the fundamental theorem
of algebra does not extend to the multi-dimensional case. It is obvious
that a necessary condition for the existence of $A^{-1}$ is $M <r$.
However this condition is no longer sufficient, since the sampling points
need not be appropriately distributed. In higher 
dimensions, the zero set of a polynomial is an algebraic curve
or an algebraic surface. For $A$ to be invertible, the samples must
not be contained in any algebraic surface. It is an open problem to
efficiently characterize all sampling sets that yield an invertible
matrix $A$.

It is still possible to obtain conditions that guarantee the existence of 
$A^{-1}$ as well as to derive estimates for the condition number of $A$ in 
the multi-dimensional case. This can be done for instance by adapting the 
approach in Section~4.3 of~\cite{Gro92} to our situation. However the 
estimates are no longer sharp and get worse with increasing dimension. We 
do not pursue this direction here.

%\section{Numerical simulations} \label{s:num}

\section{Numerical experiments: An example from geophysics} \label{s:five}

We demonstrate the performance of the proposed algorithm by applying it
to a scattered data problem from geophysics. Exploration geophysics relies on 
measurements of the Earth's physical properties like the magnetic or
gravitational field, with the goal of detecting 
anomalies which reveal underlying geological features.
In geophysical practice, it is essentially impossible to gather data in
a form that allows direct interpretation. Geoscientists, used to look
at their measurements on maps or profiles and aim at further
processing, need a representation of the originally
irregularly spaced (scattered) data points on a regular grid. 
The reconstruction or approximation of potential fields on regular grids 
from scattered data is thus one of the first and crucial steps in
the analysis of geophysical data.

As test example we use a synthetic anomaly $f$ that
represents the gravitational acceleration caused by an ensemble of buried
rectangular boxes of different size, depth, and density contrast,
see~Fig.~\ref{fig:geo1}(a). This example has also been used in~\cite{RS98}.
We sample this function at 496 randomly spaced points $(x_j,y_j)$ in the 
interval $[0,1] \times [0,1]$. Since in practice measurements are always 
contaminated by noise we add white Gaussian noise in the amount of 5\% of 
the $\ltsp$-norm of the samples $f(x_j,y_j)$. We want to reconstruct the 
function on a regular grid $\Gamma$ consisting of the grid points 
$\{(k/150, l/150)\}_{k,l =0}^{150}$.

In order to demonstrate the advantage of using Neumann boundary conditions
over periodic boundary conditions we compare the proposed algorithm to 
the so-called ACT method~\cite{FGS95,GS01}.
The latter has become a main ingredient for several 
approximation methods in geophysics~\cite{RS98,DSH99}.
We also include in the comparison 
the approximation obtained by cubic spline interpolation, which we computed
via the MATLAB function {\tt griddata} using the option 'cubic'. 

For the two methods using trigonometric approximation we use the same 
number of coefficients for the
approximating polynomial. We use a total of 11 coefficients in the
x-coordinate and the same number in the y-coordinate, resulting in
approximating polynomials of degree 121 for both methods.

Since we know the original anomaly $f$ we can compute the error between
the approximation $f_a$ and $f$ via 
$e(f_a)=\|f(\Gamma)-f_a(\Gamma)\|_2/\|f(\Gamma)\|$ on the grid $\Gamma$. 
The proposed method gives an error of $0.029$, the 
ACT method yields approximation error $0.072$, and the approximation
computed via cubic splines returns an error of $0.045$.
The approximation computed by the proposed method is appealing
both from a visual and from an approximation error viewpoint.

\begin{figure}
\begin{center}
    \subfigure[Synthetic gravity anomaly, gravity is in mGal.]{
    \epsfig{file=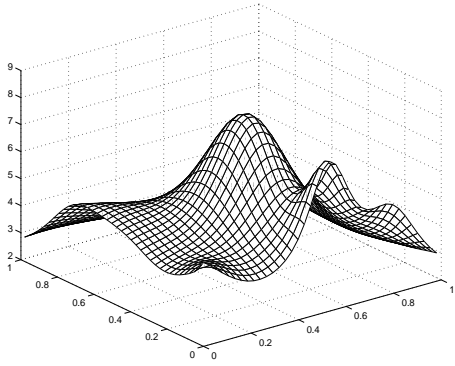,width=60mm,height=48mm}} \quad
    \subfigure[Contour plot of gravity anomaly; sampling locations
               are marked as ``o''.]{
    \epsfig{file=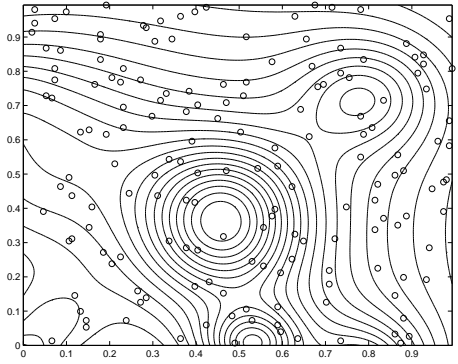,width=60mm,height=48mm}} \\
    \subfigure[ACT method, error = 0.072]{
    \epsfig{file=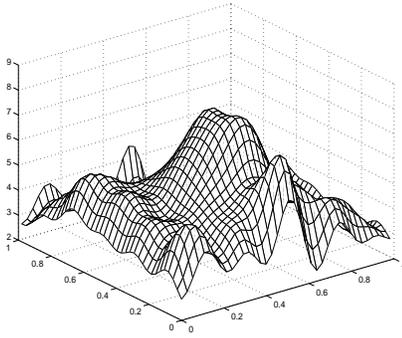,width=60mm,height=48mm}} \quad
    \subfigure[Cubic splines, error = 0.045]{
    \epsfig{file=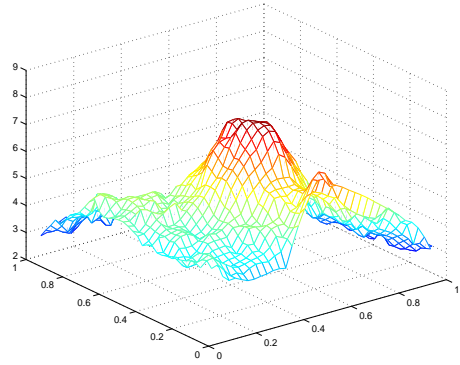,width=60mm,height=48mm}}
    \subfigure[Proposed method, error = 0.029 ]{
    \epsfig{file=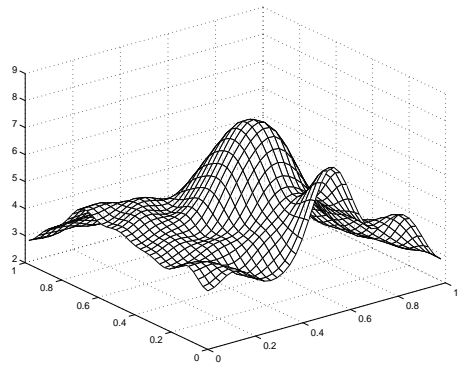,width=60mm,height=48mm}} 
    \label{fig:geo1}
    \caption{Approximation of gravitational anomaly from noisy scattered
    data (5\% noise) by proposed method and comparison to standard algorithms.}
\end{center}
\end{figure}

The significantly larger error of ACT is only due to boundary effects.
We note that there are several ways to improve the performance of the
ACT method, see~\cite{RS98}, which makes it indeed a powerful approximation 
method in geophysics~\cite{RS98,DSH99}. Since all these modifications 
can also be applied to the proposed method we expect that the proposed
(modified) algorithm will still be significantly better than
the modified ACT method.

The results of this experiment do not mean that the proposed 
method always performs better than the other two methods. Furthermore, a 
detailed comparison of various scattered data approximation methods would
have to include other standard methods such as approximation
by radial basis functions. Such a comparison is beyond the
scope of this paper.

\section*{Acknowledgment} 

T.S.\ wants to thank Raymond Chan and Michael Ng for their kind invitation
to Hongkong in December 2000. The initial steps to this research
were carried out during this wonderful and inspiring visit.

%\bibliographystyle{siam}
%\bibliography{mathbook,nuhag,linalg,samp,thomas}

\begin{thebibliography}{10}

\bibitem{BT97}
{\sc G.~Baszenski and M.~Tasche}, {\em Fast polynomial multiplication and
  convolutions related to the discrete cosine transform}, Linear Algebra Appl.,
  252 (1997), pp.~1--25.

\bibitem{DSH99}
{\sc A.~Duijndam, M.~Schonewille, and C.~Hindriks}, {\em Reconstruction of
  band-limited data irregularly sampled along one spatial direction},
  Geophysics, 64 (1999), pp.~524--538.

\bibitem{Fas97}
{\sc H.~Fassbender}, {\em On numerical methods for discrete least-squares
  approximation by trigonometric polynomials}, Math.\ Comp., 66 (1997),
  pp.~719--741.

\bibitem{FGS95}
{\sc H.~G. Feichtinger, K.~Gr{\"o}chenig, and T.~Strohmer}, {\em Efficient
  numerical methods in non-uniform sampling theory}, Numerische Mathematik, 69
  (1995), pp.~423--440.

\bibitem{GL96}
{\sc G.~Golub and C.~van Loan}, {\em Matrix Computations}, Johns Hopkins,
  Baltimore, third~ed., 1996.

\bibitem{Gro92}
{\sc K.~Gr{\"o}chenig}, {\em Reconstruction algorithms in irregular sampling},
  Math. Comp., 59 (1992), pp.~181--194.

\bibitem{Gro93a}
{\sc K.~{Gr\"ochenig}}, {\em A discrete theory of irregular sampling}, Lin.
  Alg. and Appl., 193 (1993), pp.~129--150.

\bibitem{GS01}
{\sc K.~Gr{\"o}chenig and T.~Strohmer}, {\em Numerical and theoretical aspects
  of non-uniform sampling of band-limited images}, in Theory and Practice of
  Nonuniform Sampling, F.~Marvasti, ed., Kluwer/Plenum, 2001.

\bibitem{HLP52}
{\sc G.~H. Hardy, J.~E. Littlewood, and G.~P{\'o}lya}, {\em Inequalities},
  Cambridge University Press, Cambridge, 1952.

\bibitem{Hei01}
{\sc G.~Heinig}, {\em Chebyshev-{H}ankel matrices and the splitting approach
  for centrosymmetric {T}oeplitz-plus-{H}ankel matrices}, Linear Algebra Appl.,
  327 (2001), pp.~181--196.

\bibitem{HR98}
{\sc G.~Heinig and K.~Rost}, {\em Representation of {T}oeplitz-plus-{H}ankel
  matrices using trigonometric transformations with applications to fast
  matrix-vector multiplication}, Linear Algebra Appl., 275--276 (1998),
  pp.~225--248.

\bibitem{KS99}
{\sc T.~Kailath and A.~Sayed}, {\em Fast Reliable Algorithms for Matrices with
  Structure}, SIAM, Philadelphia. PA, 1999.

\bibitem{NCT99}
{\sc M.~K. Ng, R.~H. Chan, and W.-C. Tang}, {\em A fast algorithm for
  deblurring models with {N}eumann boundary conditions}, SIAM J. Sci. Comput.,
  21 (1999), pp.~851--866 (electronic).

\bibitem{Pot01}
{\sc D.~Potts}, {\em Fast algorithms for discrete polynomial transforms on
  arbirtrary grids}, 2001.
\newblock preprint.

\bibitem{PST01}
{\sc D.~Potts, G.~Steidl, and M.~Tasche}, {\em Fast {F}ourier transforms for
  nonequispaced data: a tutorial}, in {M}odern {S}ampling {T}heory:
  {M}athematics and {A}pplications, J.~Benedetto and P.~Ferreira, eds.,
  Birkh{\"a}user, 2001, pp.~247--270.

\bibitem{RS98}
{\sc M.~Rauth and T.~Strohmer}, {\em Smooth approximation of potential fields
  from noisy scattered data}, Geophysics, 63 (1998), pp.~85--94.

\bibitem{RAG91}
{\sc L.~Reichel, G.~Ammar, and W.~Gragg}, {\em Discrete least squares
  approximation by trigonometric polynomials}, Math. Comp., 57 (1991),
  pp.~273--289.

\bibitem{SPP95}
{\sc V.~Sanchez, P.~Garcia, A.~Peinado, J.~Segura, and A.~Rubio}, {\em
  Diagonalizing properties of the discrete cosine transform}, IEEE Trans.\
  Sig.\ Proc., 43 (1995), pp.~2631--2641.

\bibitem{Sch95}
{\sc R.~Schaback}, {\em Multivariate interpolation and approximation by
  translates of a basis function}, in Approximation theory VIII, Vol.\ 1
  (College Station, TX, 1995), World Sci. Publishing, River Edge, NJ, 1995,
  pp.~491--514.

\bibitem{SS98}
{\sc O.~Scherzer and T.~Strohmer}, {\em A multi--level algorithm for the
  solution of moment problems}, Num.Funct.Anal.Opt., 19 (1998), pp.~353--375.

\bibitem{Ste92}
{\sc G.~Steidl}, {\em Fast radix-$p$ discrete cosine transform}, Appl. Algebra
  Engrg. Comm. Comput., 3 (1992), pp.~39--46.

\bibitem{Str99d}
{\sc G.~Strang}, {\em The discrete cosine transform}, SIAM Rev., 41 (1999),
  pp.~135--147 (electronic).

\bibitem{Str98}
{\sc T.~Strohmer}, {\em A {L}evinson-{G}alerkin algorithm for regularized
  trigonometric approximation}, SIAM J.\ Sci.\ Comp., 22 (2000),
  pp.~1160--1183.

\bibitem{Str99c}
\leavevmode\vrule height 2pt depth -1.6pt width 23pt, {\em Numerical analysis
  of the non-uniform sampling problem}, J.\ Comp.\ Appl.\ Math., 122 (2000),
  pp.~297--316.

\bibitem{Loa92}
{\sc C.~van Loan}, {\em Computational Frameworks for the {F}ast {F}ourier
  {T}ransform}, SIAM, 1992.

\end{thebibliography}

\end{document}